\documentclass[12pt,reqno]{amsart}
\usepackage{hyperref}
\usepackage{enumerate,fullpage}
\usepackage{amssymb,amsmath,graphicx,amsthm}
\usepackage{color}
\usepackage{soul}

\def\B{\4B}

\def\Om{\Omega}

\def\no#1{\|#1\|}

\numberwithin{equation}{section}
\def\db{\bar\partial}
\def\db*{\bar\partial^*}
\def\T{\text}

\def\simleq{\lesssim}
\def\1#1{\overline{#1}}
\def\2#1{\widetilde{#1}}
\def\3#1{\widehat{#1}}
\def\4#1{\mathbb{#1}}

\def\5#1{\frak{#1}}
\def\6#1{{\mathcal{#1}}}
\def\C{{\4C}}

\def\B{\4B}

\def\di{\partial}
\def\dib{\bar\partial}

\def\phi{\varphi}
\def\eps{\varepsilon}

\def\Om{\Omega}

\newtheorem{Thm}{Theorem}[section]
\newtheorem{Cor}[Thm]{Corollary}
\newtheorem{Pro}[Thm]{Proposition}
\newtheorem{Lem}[Thm]{Lemma}
\theoremstyle{definition}\newtheorem{Def}[Thm]{Definition}
\theoremstyle{remark}
\newtheorem{Rem}[Thm]{Remark}
\newtheorem{Exa}[Thm]{Example}

\def\Label#1{\label{#1}}
\def\bl{\begin{Lem}}
\def\el{\end{Lem}}
\def\bp{\begin{Pro}}
\def\ep{\end{Pro}}
\def\bt{\begin{Thm}}
\def\et{\end{Thm}}
\def\bc{\begin{Cor}}
\def\ec{\end{Cor}}
\def\bd{\begin{Def}}
\def\ed{\end{Def}}
\def\br{\begin{Rem}}
\def\er{\end{Rem}}
\def\be{\begin{Exa}}
\def\ee{\end{Exa}}
\def\bpf{\begin{proof}}
\def\epf{\end{proof}}
\def\ben{\begin{enumerate}}
\def\een{\end{enumerate}}

\def\1alpha{[\frac1\alpha]}
\def\T{\text}

\def\C{{\mathbb C}}

\begin{document}
\title[The complex Monge-Amp\`ere...]{The complex Monge-Amp\`ere equation on weakly pseudoconvex domains}
\author[L.~Baracco, T.V.~Khanh, S.~Pinton]{Luca Baracco, Tran Vu Khanh and Stefano Pinton}
\address {Luca Baracco, Stefano Pinton}
\address{Dipartimento di Matematica, Universit\`a di Padova, via 
Trieste 63, 35121 Padova, Italy}
\email{baracco@math.unipd.it,pinton@math.unipd.it}

	\address{Tran Vu Khanh}
	\address{School of Mathematics and Applied Statistics, University of Wollongong, NSW, Australia,  2522}
	\email{tkhanh@uow.edu.au}
\thanks{The research of T.~V.~Khanh was supported by the Australian Research Council
	DE160100173.}

\abstract
We show  here  a ``weak" H\"older-regularity up to the boundary of the solution to the Dirichlet problem for the complex Monge-Amp\`{e}re equation with data in the $L^p$ space and $\Om$ satisfying an $f$-property.  The $f$-property is a potential-theoretical condition which holds for all pseudoconvex domains of finite type and many examples of infinite type. 
\endabstract
\maketitle
%\endtopmatter

\section{Introduction}
For a $C^2$, bounded, pseudoconvex domain $\Om\subset\subset \C^n$, the Dirichlet problem for the Monge-Amp\`ere equation consists in
\begin{equation}
\Label{1.1}
\begin{cases}
u\in PSH(\Omega)\cap L^\infty_{loc}(\Omega),
\\
(dd^cu)^n=\psi dV&\T{in }\Om,
\\
u=\phi&\T{on }b\Om.
\end{cases}
\end{equation}
 A great deal of work has been done when $\Om$ is strongly pseudoconvex. 
  On this domain, we can divide the literature by three kinds of data $\psi$.
\begin{itemize}
  \item \underline{The H\"older data:} Bedford-Taylor prove in \cite{BeTa76} that $u\in C^{\frac{\alpha}{2}}(\bar\Om)$ if $\phi\in C^{\alpha}(b\Om)$, $\psi^{\frac1n}\in C^{\frac\alpha2}(\bar\Om)$ for $0<\alpha\le 2$.
  \item \underline{The smooth data:} Caffarelli, Kohn and Nirenberg prove in \cite{CKNS85} that $u\in C^\infty(\bar{\Om})$, for $\phi\in C^\infty(b\Om)$ and $\psi\in C^\infty(\bar{\Om})$, in case $\psi>0$ in $\bar\Om$ and $b\Om$ is smooth. 
  \item \underline{The $L^p$ data:} Guedj, Kolodziej and Zeriahi prove in \cite{GKZ08} that if $\psi\in L^p(\Om)$ with $p>1$ and $\phi\in C^{1,1}(b\Om)$ then $u\in C^{\gamma}(\bar\Om)$ for any $\gamma<\gamma_p:=\frac2{qn+1}$ where $\frac1q+\frac1p=1$.
\end{itemize}  
 When $\Om$ is no longer strongly pseudoconvex but has a certain ``finite type", there are some known results for this problem due to Blocki \cite{Blo96}, Coman \cite{Com97}, and Li \cite{Li04}. Recently, Ha and the second author gave a general related result to a H\"older data under the hypothesis that $\Om$ satisfies an $f$-property (see Definition~\ref{d21} below).  The $f$-property is a consequence of the geometric ``type"  of the boundary. All pseudoconvex domains of finite type satisfy the $f$-property as well as 
 many  classes of domains of infinite type (see \cite{HaKh15, Kha16a, KhZa10} for discussion on the $f$-property).  Using the $f$-property, a  ``weak" H\"older regularity for the solution of the Dirichlet problem of complex Monge-Amp\`ere equation is obtained in \cite{HaKh15}.  
Coming back to the case of $\Om$ of finite type, in a recent paper with Zampieri \cite{BaKhPiZa16}, we prove the H\"older regularity for $\psi\in L^p$ with $p>1$. The purpose of the present paper is to generalize the result in \cite{BaKhPiZa16} to a pseudoconvex domain satisfying an $f$-property.  For this purpose we recall the definition of a weak H\"older space in \cite{ HaKh15,Kha16a}. 
Let $f$ be an increasing function such that $\underset{t\to+\infty}{\lim} f(t)=+\infty$, $f(t)\lesssim t$. For a subset $A$ of $\C^n$, define the $f$-H\"older space on $A$ by
$$\Lambda^{f}(A)=\{u : \no{u}_{L^\infty(A)}+\sup_{z,w\in A, z\not=w}f(|z-w|^{-1}) \cdot |u(z)-u(w)|<\infty \}$$ 
and set 
$$\no{u}_{\Lambda^f(A)}= \no{u}_{L^\infty(A)}+\sup_{z,w\in A, z\not=w}f(|z-w|^{-1})\cdot |u(z)-u(w)|. $$
Note that the notion of the $f$-H\"older space includes the standard H\"older space $\Lambda_\alpha$ by taking $f(t) = t^{\alpha}$ (so that $f(|h|^{-1}) = |h|^{-\alpha}$) with 
$0<\alpha\le1$. 
Here is our result
\begin{Thm}\label{t1} Let $\Omega\subset\C^n$ be a bounded, pseudoconvex domain admitting the $f$-property. Suppose that $
 \displaystyle\int_1^\infty \dfrac{da}{a f(a)}<\infty$ and denote by $g(t):=\left(\displaystyle\int_t^\infty \dfrac{da}{a f(a)}\right)^{-1}$ for $t\ge 1$. 
If $0<\alpha\le 2$, $\phi\in \Lambda^{t^{\alpha}}(b\Omega)$, 
and $\psi\ge 0$ on $\Omega$ with $\psi\in L^p$ with $p>1$, then the  Dirichlet problem for the complex Monge-Amp\`ere equation \eqref{1.1}
has a unique plurisubharmonic solution $u\in \Lambda^{g^{\beta}}(\overline{\Omega})$. Here $\beta=\min(\alpha,\gamma)$, for any $\gamma<\gamma_p=\frac{2}{nq+1}$ where $\frac1p+\frac1q=1$.
\end{Thm}
The proof follows immediately from Theorem~\ref{defining} and \ref{t23} below.
Throughout the paper we use $\lesssim$ and $\gtrsim$ to denote an estimate up to a positive constant, and $\approx$ when both of them hold
simultaneously. Finally, the indices $p$, $\alpha$, $\beta$, $\gamma$ and $\gamma_p$ only take ranges as in Theorem~\ref{t1}.  

\section{H\"older regularity of the solution}
\Label{s2}
We start this section by defining the $f$-property as in \cite{Kha16a,KhZa10}.
\begin{Def}\label{d21}
\label{d1}  For a smooth, monotonic,
increasing function $f :[1,+\infty)\to[1,+\infty)$ with ${f (t)}{t^{-1/2}}$ decreasing, we say that $\Om$ has the $f$-property if there exist a neighborhood $U$ of $b\Om$ and a family of functions $\{\phi_\delta\}$ such that
\begin{enumerate}
  \item [(i)] the functions $\phi_\delta$ are plurisubharmonic,  $C^2$ on $U$, and satisfy $-1\le \phi_\delta \le0$,
  \item[(ii)] $i\di\dib \phi_\delta\gtrsim f(\delta^{-1})^2Id$ and $|D\phi_\delta|\lesssim  \delta^{-1}$ for any  $z\in U\cap \{z\in \Om:-\delta<r(z)<0\}$, where $r$ is a $C^2$-defining function of $\Om$.
  \end{enumerate}
In \cite{Kha16a}, using the $f$-property, the second author constructed a family of plurisubharmonic peak functions with good estimates. This family of plurisubhamonic peak functions yields the existence of a defining function $\rho$ which is uniformly strictly-plurisubharmonic and weakly H\"older (see \cite{HaKh15}).
\end{Def} 
\begin{Thm}[Khanh \cite{Kha16a} and Ha-Khanh\cite{HaKh15}]\Label{defining} Assume that $\Omega$ is a bounded, pseudoconvex domain  admitting the $f$-property as in Theorem~\ref{t1}. Then there exists a uniformly strictly-plurisubharmonic defining function of $\Om$ which belongs to the $g^2$-H\"older space of $\overline{\Om}$, that means, 
\begin{eqnarray}\label{rho}
\rho\in\Lambda^{g^2}(\bar\Om),\quad\Om=\{\rho<0\}  \quad\T{and}\quad i\di\dib \rho\ge \T{Id}.
\end{eqnarray}
\end{Thm}
The existence and uniqueness of the solution $u\in L^\infty(\Om)$ to the equation \eqref{1.1} need a weaker condition, in particular, one only need $\rho\in C^0(\bar\Om)$ as shown by \cite{Kolo98}.  
\bt[Kolodziej \cite{Kolo98}]\Label{t3.0} Let $\Om$ be a bounded domain in $\C^n$. Assume that there exists a function $\rho$ such that 
$$\rho\in C^0(\bar\Om), \quad \Om=\{\rho<0\}\quad\T{and}\quad i\di\dib\rho\ge \T{Id}.$$
 Then for any $\phi\in C^0(b\Om)$, $\psi\in L^p(\Om)$ there is a unique plurisubharmonic solution $u(\Om,\phi,\psi)\in C^0(\bar\Om)$.
\et
To improve the smoothness of $u$, we increase the smoothness of $\rho$ and $\psi$. \begin{Thm}[Ha-Khanh \cite{HaKh15}]\label{t22} Let $\rho$ satisfy \eqref{rho}. If $\phi\in \Lambda^{t^{\alpha}}(b\Omega)$ and  $\psi^{\frac1n}\in \Lambda^{g^{\alpha}}(\overline{\Omega})$, then the  Dirichlet problem for the complex Monge-Amp\`ere equation \eqref{1.1}
has a unique plurisubharmonic solution $u(\Om,\phi,\psi)\in \Lambda^{g^{\alpha}}(\overline{\Omega})$.
\end{Thm}
Now we focus on lowering the smoothness of $\psi$ and prove the following theorem.
\begin{Thm}\Label{t23}
 Let $\rho$ satisfy \eqref{rho}. If  $\phi\in \Lambda^{t^{\alpha}}(b\Omega)$ and  $\psi\in L^p(\Omega)$, then the  Dirichlet problem for the complex Monge-Amp\`ere equation \eqref{1.1}
has a unique plurisubharmonic solution $u(\Om,\phi,\psi)\in \Lambda^{g^{\beta}}(\overline{\Omega})$.
\end{Thm}
In order to prove this theorem, we need to construct a subsolution  with $L^p$ data. Here, $v$ is a subsolution to \eqref{1.1} in the sense that $v$ is plurisubharmonic, $v|_{b\Om}=\phi$ and $(dd^cv)^n\ge \psi dV$ in $\Om$. 
\bp
\Label{p2.3}
Let $\rho$ satisfy \eqref{rho}.Then there is a subsolution $v\in \Lambda^{g^\beta}(\bar\Om)$ to \eqref{1.1} for $\phi\in C^\alpha(b{\Om})$ and $\psi\in L^p(\Om)$.
\ep

\bpf
For a large ball $\B$ containing $\Om$,  we set  $
\tilde \psi(z):=\begin{cases} \psi(z)\quad\T{if $z\in\Om$},
\\
0\quad\T{ if $z\in\B\setminus\Om$}.
\end{cases}
$
 First, we apply  Theorem 1 in \cite{GKZ08} on $\B$ with $\tilde \psi\in L^p(\B)$ and zero-valued boundary condition, it follows $u_1=u(\B,0,\tilde \psi)\in \Lambda^{t^\gamma}(\bar{\B})$. Second, we apply Theorem~\ref{t22} on $\Om$ twice: first for $u_2:=u(\Om,-u_1|_{b\Om}, 0)\in  \Lambda^{g^\gamma}$ since $u_1|_{b\Om}\in \Lambda^{t^\gamma}$  and second for $u_3:=u(\Om,\phi,0)\in \Lambda^{g^\alpha}$ by the hypothesis $\phi\in \Lambda^{t^\alpha}$. Finally,  taking the summation $v=u_1+u_2+u_3$, we  have the conclusion.

\epf

\bpf[Proof of Theorem~\ref{t23}]
Keeping the notation of Theorem~\ref{t3.0}, let $u(\Om,\phi,\psi)\in C^0(\bar\Om)$ be the solution of \eqref{1.1}. What follows  is dedicated to showing that this $C^0$ plurisubharmonic solution $u(\Om,\phi,\psi)$ is in fact in $\Lambda^{g^\beta}(\bar\Om)$. By Theorem~\ref{t22}, $w:=u(\Omega,\phi,0)\in \Lambda^{g^\alpha}(\bar\Om)$; let $v$ be as in Proposition~\ref{p2.3}, comparison principle yields at once
\begin{equation}
\Label{3.1}
v\le u(\Omega,\phi,\psi)\le w.
\end{equation}
By \eqref{3.1} and the $g^\beta$-H\"older regularity of $v$ and $w$ we get
$$
|u(z)-u(\zeta)|\simleq [g(|z-\zeta|^{-1})]^{-\beta}\quad z\in \bar\Om,\,\,\zeta\in b\Om,
$$
and therefore for $\delta$ suitably small
\begin{equation}
\Label{3.0}
|u(z)-u(z')|\simleq [g(\delta^{-1})]^{-\beta},\quad z,\,z'\in\Om\setminus \Om_\delta\T{ and }|z-z'|<\delta
\end{equation}
 where $\Omega_{\delta}:=\{ z\in\C^n:\, r(z)<-\delta \}$ and $r$ is the $C^2$ defining function for $\Omega$ with $|\nabla r|=1$  on $b\Omega$. We have to prove that \eqref{3.0} also holds for $z,\,z'\in\Om_\delta$. 
For $z\in\bar\Om_\delta$, we use the notation
\begin{equation*}
u_{\frac\delta2}(z):=\underset{|\zeta|<{\frac\delta2}}\sup u(z+\zeta),\qquad
\tilde u_{\frac\delta2}(z):=\frac1{\sigma_{2n-1}{\big(\frac\delta2}\big)^{2n-1}}\int_{b\B(z,{\frac\delta2})}u(\zeta)dS(\zeta),
\end{equation*}
and 
$$
\hat u_{\frac\delta2}(z):=\frac1{\sigma_{2n}{(\frac\delta2})^{2n}}\int_{\B(z,{\frac\delta2})}u(\zeta)dV(\zeta),
$$
where ${\sigma_{2n-1}{\big(\frac\delta2}\big)^{2n-1}}=Vol(b\B(z,\frac\delta2))$ and ${\sigma_{2n}{\big(\frac\delta2}\big)^{2n}}=Vol(\B(z,\frac\delta2))$. It is  obvious that
\begin{eqnarray}\label{order}
\hat u_{\frac\delta2}\le \tilde u_{\frac\delta2}\le u_{\frac\delta2} \quad\T{in}\quad\Om_\delta.
\end{eqnarray}
Furthermore, we have $L^1$ estimate of the difference between $u$ and $\tilde{u}_{\frac{\delta}{2}}$ and the stability estimate in the following theorems:

\bt[Baracco-Khanh-Pinton-Zampieri \cite{BaKhPiZa16}]
\Label{t3.1}
For any $0<\epsilon<1$, we have
\begin{equation}
\Label{3.2}
\no{\tilde u_{\frac\delta2}-u}_{L^1(\Om_\delta)}\simleq \delta^{1-\epsilon}.
\end{equation}
\et

\bt[Guedj-Kolodziej-Zeriahi \cite{GKZ08}]\label{tse}
Fix $0\leq f\in L^p(\Omega),\, p>1$. Let $U, \, W$ be two bounded plurisubharmonic functions in $\Omega$ such that $(dd^c U)^n=fdV$ in $\Omega$ and let $U\geq W$ on $\partial\Omega$. Fix $s\geq 1$ and $0\leq \eta< \frac{s}{nq+s}$, $\frac 1p+\frac 1q=1$. Then there exists a uniform constant $C=C(\eta, \|f\|_{L^p(\Omega)})>0$ such that 
$$
\sup_\Om(W-U)\leq C\no{(W-U)_+}^\eta_{L^s(\Om)},
$$
where $(W-U)_+:=\max(W-U,0)$.
\et

By \eqref{3.0}, we have 
$$\tilde u_{\frac\delta2}\le u_{\frac\delta2}\le u+ c[g(\delta^{-1})]^{-\beta}, \quad \T{on $b\Om_\delta$ for suitable constant $c$}.$$ 
Thus, we can  apply Theorem \ref{tse} for $\Om_\delta$ with $U:=u+c[g(\delta^{-1})]^{-\beta}$, $W: =\tilde u_{\frac\delta2}$ and $s:=1$; thus we get

\begin{equation}
\Label{3.7}
\begin{split}
\underset{\Om_\delta}\sup\left(\tilde u_{\frac\delta2}-(u+c[g(\delta^{-1})]^{-\beta})\right)&\underset{\T{Theorem~\ref{tse}}}\lesssim \no{\left(\tilde u_{\frac\delta2}-(u+c[g(\delta^{-1})]^{-\beta})\right )_+}^\eta_{L^1(\Om_\delta)}
\\
&\lesssim \no{\tilde u_{\frac\delta2}-u}^\eta_{L^1(\Om_\delta)}\underset{\T{Theorem~\ref{t3.1}}}\lesssim \delta^{(1-\epsilon)\eta},
\end{split}
\end{equation}
for any $\eta<\frac{1}{2}\gamma_p=\frac{1}{nq+1}$ where $\frac 1q+\frac 1p=1$.
Taking $\gamma<\gamma_p$, $\beta=\min(\alpha,\gamma)$, $\epsilon=\frac{\gamma_p-\gamma}{\gamma_p+\gamma}>0$ and $\eta=\frac{1}{4}(\gamma_p+\gamma)< \frac{1}{2}\gamma_p$ so that $(1-\eps)\eta=\frac{\gamma}{2}$, it follows
\begin{equation}\label{a1}
\underset{\Om_\delta}\sup\left(\tilde u_{\frac\delta2}-u\right)\lesssim\delta^{(1-\epsilon)\eta}+[g(\delta^{-1})]^{-\beta}\lesssim \delta^{\frac{\gamma}{2}}+ [g(\delta^{-1})]^{-\beta}\lesssim [g(\delta^{-1})]^{-\beta},
\end{equation}
where the last inequality of \eqref{a1} follows by $g(\delta^{-1})\lesssim \delta^{-\frac12}$ (by the conditions on $f$ in the $f$-property).						\\
 
Similarly to \cite[Lemma 4.2]{GKZ08} by using the fact that  $g(c\delta^{-1}))\approx g(\delta^{-1})$ for any constant $c>0$, one can state the equivalence between
		$$\sup_{\Om_\delta}(u_\delta-u)\lesssim  [g(\delta^{-1})]^{-\beta}\quad\T{and}\quad 
		\sup_{\Om_\delta}(\hat{u}_\delta-u)\lesssim [g(\delta^{-1})]^{-\beta}.$$
Using this equivalence together with the inequalities in \eqref{order}, it follows that \eqref{a1} is equivalent to 
\begin{equation}
\Label{b}
\underset {\Om_\delta}\sup(u_\frac\delta2-u)\lesssim [g(\delta^{-1})]^{-\beta}.
\end{equation}
From \eqref{3.0} and \eqref{b}, it is easy to prove that 
$$|u(z)-u(z')|\lesssim [g(|z-z'|^{-1})]^{-\beta}\quad\T{for any $z,z'\in \bar{\Om}$}.$$ 
\epf

\end{document}